\documentclass{ifacconf}

\usepackage{graphicx}      
\usepackage{natbib}        
\usepackage{amsmath, amssymb}
\usepackage{mathtools}
\mathtoolsset{showonlyrefs}
\usepackage{algorithm}
\usepackage[]{algpseudocode}
\usepackage{color}
\usepackage{tikz}
\usepackage{subcaption}
\usepackage{environ}

\makeatletter
\def\BState{\State\hskip-\ALG@thistlm}
\newsavebox{\measure@tikzpicture}
\NewEnviron{scaletikzpicturetowidth}[1]{%
  \def\tikz@width{#1}%
  \def\tikzscale{1}\begin{lrbox}{\measure@tikzpicture}%
  \BODY
  \end{lrbox}%
  \pgfmathparse{#1/\wd\measure@tikzpicture}%
  \edef\tikzscale{\pgfmathresult}%
  \BODY
}
\makeatother
\usepackage{xcolor}
\definecolor{edgecl}{HTML}{E27D60}
\definecolor{treecl}{HTML}{41B3A3}
\definecolor{nodecl}{HTML}{E8A87C}

\newcommand{\lap}[0]{{\mathcal{L}}}

\newcommand{\graph}[0]{{\mathcal{G}}}

\newcommand{\Tr}{\mathbf{tr}}

\newcommand{\Verts}{\mathcal{V}}

\newcommand{\Edges}{\mathcal{E}}
\newcommand{\Weights}{\mathcal{W}}
\newcommand{\tree}{{\mathcal{T}}}
\newcommand{\cycle}{{\mathcal{C}}}

\newcommand{\Htwo}{{\mathcal{H}_2}}
\newcommand{\Htwosq}{{\mathcal{H}_2^2}}

\begin{document}

\def\q{\quad}

\long\def\okkvir#1#2#3{%
\hbox{\vrule\vbox{\hrule\vskip#1pt\hbox{\hskip#1pt\vbox{%
\hsize#2cm\parindent=0pt #3%
}\hskip#1pt}\vskip#1pt\hrule}\vrule}%
}

\begin{frontmatter}

\title{Graph-Theoretic Optimization\\ for  Edge Consensus} 

\author[First]{Mathias Hudoba de Badyn}
\author[Second]{Dillon R. Foight} 
\author[Second]{Daniel Calderone}
\author[Second]{Mehran Mesbahi}
\author[First]{Roy S. Smith}

\address[First]{Automatic Control Laboratory, ETH Z\"{u}rich, Switzerland\\ (e-mails: { \{mbadyn,~rsmith\}@control.ee.ethz.ch})}
\address[Second]{William E. Boeing Department of Aeronautics \& Astronautics, University of Washington, Seattle WA, 98109 USA\\ (e-mails: {\{dfoight,~djcal,~mesbahi\}@uw.edu})}

\begin{abstract}
  We consider network structures that optimize the $\Htwo$-norm of weighted, time scaled consensus networks, under a minimal representation of such consensus networks described by the edge Laplacian.
  We show that a greedy algorithm can be used to find the minimum-$\Htwo$ norm spanning tree, as well as how to choose edges to optimize the $\Htwo$ norm when edges are added back to a spanning tree.
  In the case of edge consensus with a measurement model considering all edges in the graph, we show that adding edges between slow nodes in the graph provides the smallest increase in the $\Htwo$ norm. 
\end{abstract}


\end{frontmatter}







\section{Introduction}

Many natural and synthetic systems are distributed among agents in a network.
One popular information-sharing protocol distributed over a network is \emph{consensus}, used in many applications ranging from robotics, as in~\cite{Joordens2009}; sensor networks, as in~\cite{Olfati-Saber2005}; and multi-agent systems, as in~\cite{Olfati-Saber2004}, and~\cite{Tanner2004}.
Much work on studying networked dynamical systems has focused on how the physical structure of the network affects its dynamics.
For example, network symmetries have been shown to be related to controllability of consensus, as discussed by~\cite{Rahmani2009a} and \cite{Alemzadeh2017}.

In this work, we are interested in studying system-theoretic measures of performance, and how to optimize graph structures for these system-theoretic measures.
In particular, we examine the $\mathcal{H}_2$ norm, which can be interpreted as a measure of noise attenuation over the network, as discussed by \cite{Siami2014a}, and has a rich history of use in the literature.
For consensus networks with leaders, the $\Htwo$ norm has an interpretation in terms of effective resistance across the network, as discussed in~\cite{Chapman2013a}, and in~\cite{Hudobadebadyn2019,HudobadeBadyn2021}, a fast method of computing the $\Htwo$ norm for series-parallel networks exploited this interpretation. 
\cite{Bamieh2012} and \cite{Patterson2010a,Patterson2014} have utilized the related concept of coherence to consider local feedback laws and leader selection to promote coherence.

A minimal representation of consensus networks can be obtained by looking at a system where only relative agent states across connections are considered.
Such a system, represented by the \emph{edge Laplacian} rather than the \emph{graph Laplacian}, was considered for estimation and control in relative sensing networks in~\cite{Sandhu2005,Sandhu2009}, and for formation flight in~\cite{Smith2005}.
When considering weighted consensus networks operating with agents that have non-heterogeneous time scales, it was shown in~\cite{Foight2019a} that edge consensus allows explicit calculations of the $\Htwo$ norm of the system.

This latter formulation is the setting of the present work.
In this paper, we focus on the optimization of graph structures in the setting of edge-weighted, and time scaled networks with the $\Htwo$ norm of edge consensus as our optimization measure.
The contributions of the paper are as follows.
First, we show that a greedy algorithm can be used to find the minimum-$\Htwo$ norm spanning tree.
This has applications to systems where resilience to noise is desired, but a minimum number of communication links is also desired, such as in stealthy UAV systems.
We then discuss how to add weighted edges back to the graph to improve the $\mathcal{H}_2$ norm.
  In the case of edge consensus with a measurement model considering all edges in the graph, we show that adding edges between slow nodes in the graph provides the smallest increase in the $\Htwo$ norm. 

The organization of this paper is as follows.
In \S\ref{sec:prel-probl-stat}, we outline our notation, the mathematical preliminaries on graph and edge Laplacians, and describe our problem statement in the context of the $\Htwo$ performance of edge consensus.
Our main results are outlined in \S\ref{sec:main-results}, and examples are shown in \S\ref{sec:examples}.
The paper is concluded in \S\ref{sec:conclusion}.


\section{Preliminaries}
\label{sec:prel-probl-stat}

In this section, we outline notation and the mathematical preliminaries on graph and systems theory, and then lay out the setting of the problem statement.

\subsection{Preliminaries on Graphs \& Consensus}

We denote the real numbers as $\mathbb{R}$, the non-negative reals as $\mathbb{R}_+$, the positive reals as $\mathbb{R}_{++}$, and the real $n$-dimensional Euclidean vector space as $\mathbb{R}^n$.
Vectors in $\mathbb{R}^n$ are written in lower-case $x,y,z$, etc., and matrices in $\mathbb{R}^{n\times m}$ are written in capital-case $M,N,R$, etc.
$I$ denotes the identity matrix.
$\mathbf{1}\in\mathbb{R}^n$ denotes the length-$n$ vector of ones, and $\mathbf{0}$ denotes a matrix of zeros of comfortable dimensions.

A \emph{graph} $\graph$ with $n$ \emph{vertices} (or nodes) is a triple of sets $(\Verts,\Edges,\Weights)$, where $\Verts = [n]$ is a set of labeled vertices, $\Edges \subseteq \Verts\times\Verts$ is a set of $m$ \emph{edges} denoting connections between vertices, and $\Weights \in \mathbb{R}^{m}_{++}$ is a set of \emph{edge weights}, denoting some notion of  `strength' of the corresponding edge.
A graph $\graph$ is \emph{undirected} if edge $ij = ji$, and \emph{directed} otherwise.
The \emph{neighbourhood} of vertex $i$ is $N_i=\{j:ij\in\Edges\}$, and the \emph{unweighted degree} of $i$ as $\mathrm{deg}(i) = |N(i)|$.

A \emph{path} $\mathcal{P}$ is a set of edges $\{i_1i_2,~i_2i_3,\dots,i_{l-2}i_{l-1},~i_{l-1}i_l\}$ where $i_1,\dots,i_l$ are distinct,  and a \emph{cycle} $\cycle$ is a path, except that the first vertex and the last vertex are the same, i.e. $i_1 = i_l$.
An \emph{connected graph} $\graph$ is one in which there is an undirected path between any two vertices in $\graph$.
A \emph{tree} $\tree$ is a connected graph with no cycles, and a \emph{spanning tree} of a connected graph $\graph$ is a tree on the same vertex set as $\graph$ (i.e., $\Verts(\tree) = \Verts(\graph)$), with edges $\Edges(\tree) \subseteq \Edges(\graph)$.
If $\tree$ is a spanning tree of $\graph$, then we can write $\tree\subseteq \graph$, where `$\subseteq$' is to be interpreted setwise.
Note that a tree $\tree$ with $n$ vertices must have exactly $n-1$ edges.

To each vertex $i\in\Verts$, we assign a state $x_i \in \mathbb{R}$.
The dynamics of each vertex state $x_i$ are assumed, unless otherwise stated, to be single-integrator dynamics $\dot{x}_i = u_i$.
By setting $u_i$ to be a weighted average of the states of its neighbours, we arrive at the \emph{consensus dynamics}:
\begin{align}
  \dot x_i = \sum_{j\in N_i} w_{ij}(x_j - x_i).\label{eq:1}
\end{align}


One can also define the \emph{incidence matrix} $D_\graph \in \mathbb{R}^{n\times m}$ of an undirected graph $\graph$, where each column $d_{ij}$ of $D_\graph$ corresponds to an edge $ij \in \Edges$.
The column $d_{ij}$ has 1 in the $i$th position, and $-1$ in the $j$th position; the choice of setting $i=1$ or $j=1$ is arbitrary.
By setting the matrix $W_\graph \in \mathbb{R}^{m\times m}$ to be the diagonal matrix containing the edge weights $w_{ij}$, one can write the (undirected) graph Laplacian as $\lap_\graph = D_\graph W_\graph D_\graph ^\intercal $.

By appropriately stacking $x_i$ into a vector $x\in \mathbb{R}^{n}$, the consensus dynamics in \eqref{eq:1} can be written in vector form as as
\begin{align}
  \dot x = - \lap_\graph x.
\end{align}

A \emph{time-scaled graph} $\graph$ is a quadruple of sets $(\Verts,\Edges,\Weights,\mathcal{S})$, where $\Verts,\Edges$, \& $\Weights$ are as before, and $\mathcal{S} \in \mathbb{R}_{++}^n$ is a set of \emph{time scales} $\{\epsilon_i\}_{i=1}^n$ on the vertices.
We can then define the \emph{time scaled consensus dynamics} as
\begin{align}
  \epsilon_i \dot{x}_i = \sum_{j\in N_i} w_{ij}(x_j - x_i) ~\longleftrightarrow ~E_{\graph}\dot{x}_i = -\lap_{\graph} x,
\end{align}
where $E_\graph = \mathrm{diag}(\{\epsilon_i\}_{i=1}^n)$.

\subsection{Edge Consensus and Problem Configuration}

The problem considered in this paper is derived from the preceding consensus dynamics, but with process and measurement noise.
Consider a graph $\graph$ of $n$ multi-scale integrators, and zero-mean Gaussian process noise $\omega_i(t)$,
\begin{align}
\epsilon_i \dot{x}_i(t) =  u_i(t) + \omega_i(t), \label{eq:consensus}
\end{align}
where $\mathbf{E}\left[\omega(t) \omega(t)^\intercal \right] = \mathrm{diag}(\sigma_{\omega_i}^2) \triangleq\Omega^2 $, $x_i$ is the (scalar) state of the $i$-th agent, $\epsilon_i$ is node $i$'s time scale parameter, and $u_i$ is the control input.
A weighted, noisy, decentralized feedback controller that seeks to bring agents to consensus is given by,
\begin{align}
u_i(t) & = \sum_{j\in N(i)} \left[w_{ij} (x_j(t) - x_i(t)) + v_{ij}(t)\right] \nonumber \\
u(t) & =  D_\graph W_\graph D_\graph^\intercal  x(t) + D_\graph v(t), \label{eq:control}
\end{align}
where the noise over the edge $ij$ is $v_{ij}(t)$, with covariance $\mathbf{E}[v(t) v(t)^\intercal ] = \mathrm{diag}(\sigma_{v_{ij}}^2) \triangleq\Gamma^2$ for all $(i,j) \in \Edges$.
 $W_\graph\triangleq\mathrm{diag}(\mathcal{W})$ is the matrix of edge weights, and  $v(t),u(t)$ are the stacked vectors of measurement noises and control inputs.
Applying~\eqref{eq:control} to the stacked-vector version of~\eqref{eq:consensus} gives a general, time scaled and weighted consensus dynamics with process and measurement noise,
\begin{align} \label{eq:model}
\Sigma_\graph \triangleq   
  \begin{cases}
    \dot{x}(t) =  E_{\graph}^{-1} \lap_{\graph} x(t) \\ ~~~~~+ \begin{bmatrix}E_{\graph}^{-1} &  E_{\graph}^{-1} D_\graph\end{bmatrix} \begin{bmatrix} \omega(t) \\ v(t) \end{bmatrix}.
  \end{cases}
\end{align}

 The scaled Laplacian $E_{\graph}^{-1} \lap_{\graph}$ is rank-deficient (it posesses a zero eigenvalue for each connected component of $\graph$), which complicates reasoning about the $\Htwo$ performance of $\Sigma_\graph$.
 Solutions to this problem in the literature include studying the $\Htwo$ norm of the Dirichlet Laplacian, as in \cite{Chapman2015,Hudobadebadyn2019}, or as considered in this paper, a minimal representation of $\Sigma_\graph$ obtained by studying edge consensus, seen in \cite{Zelazo2013,Zelazo2011a,Foight2019a}. In the case of a time scaled graph, the \emph{edge Laplacian} is defined by $\lap^\graph \triangleq D_\graph^\intercal  E_{\graph}^{-1} D_\graph$.
We denote the {edge Laplacian} of a graph $\graph$ with a superscript: $\lap^\graph$, and the {graph Laplacian} with a subscript: $\lap_\graph$.

In~\cite{Zelazo2011a}, it was shown that the graph and edge Laplacians are related by a similarity transformation;~\cite{Foight2019a} extended  the transformation to time-scaled and matrix-weighted graphs.
This similarity transformation is constructed by choosing a spanning tree $\tree\subseteq\graph$, and noting that the cycles can be considered as linear combinations of the spanning tree edges in the following sense.
Let $D_\tree$ denote the incidence matrix of the spanning tree $\tree$, and let  $D_{\graph\setminus\tree}$  denote the incidence matrix of the remaining edges that correspond to the edges completing the cycles in $\graph$.
Define the \emph{Tucker representation} matrix $R_\graph \triangleq [I~~T_\tree]$, where  
\begin{align}
    T_\tree = \left[ D_\tree^\intercal D_\tree \right]^{-1} D_\tree ^\intercal D_{\graph\setminus\tree}.
\end{align}
We refer to the matrix $T_\tree$ as the \emph{cycle representation matrix}.
Note that if the graph is a tree, i.e. $\graph = \tree$, then the only choice of spanning tree for $R_\graph$ is in fact $\tree$, and $T_\tree$ is the empty matrix ($R_\graph = I$ in this case).


It suffices to measure the relative states of the nodes $x_i-x_j$ over the edges $ij$ of a spanning tree $\tree\subseteq\graph$ to capture consensus: if $x_i-x_j=0$ for all $ij\in\tree$, then we know that the network is at consensus.
This is the motivation for the edge consensus model, formalized below.
We define the \emph{edge states} for a edge weighted and time scaled graph as $x_e \triangleq S_v^{-1} x$, where the similarity transformation and correponding state matrix are,
\begin{align}
    S_v&= \begin{bmatrix}
             E_{\graph}^{-1} D_\tree \left[D_\tree^\intercal  E_{\graph}^{-1} D_\tree\right]^{-1} &~ \mathbf{1}
         \end{bmatrix}\\
    S_v^{-1} \lap_\graph S_v &= \begin{bmatrix}
                                   \lap^\tree R_\graph W_\graph R_\graph^\intercal &0\\0&0
                               \end{bmatrix},
\end{align}
where $\lap^\tree = D_\tree^\intercal E^{-1}_\graph D_\tree$.

The dynamics of the edge states $x_e$ corresponding to \eqref{eq:model} are then given by
\begin{align}
  \dot x_e =& \begin{bmatrix}
                 -\lap^{\tree} R_\graph W_\graph R_\graph^\intercal  & 0 \\
                 0 & 0
              \end{bmatrix}x_e 
             + \begin{bmatrix}
                    D_\tree^\intercal  E_{\graph}^{-1} &  - \lap^\tree R_\graph\\
                    \mathbf{tr}[E_\graph]^{-1}\mathbf{1}^\intercal  & 0
                \end{bmatrix}
             \begin{bmatrix}
                 \omega(t) \\ v(t)
             \end{bmatrix}.       \label{eq:16}
\end{align}
The form of the system matrix in \eqref{eq:16} suggests a partitioning of $x_e$ into the edge states of the chosen spanning tree $\tree$ and, and the consensus subspace, i.e. $x_e = [x_\tree^\intercal  ~~ x_{\mathbf{1}}]^\intercal $.
Finally, we arrive at two models of the spanning tree edge states, differing in their choice of output:
\begin{align}
    \Sigma^\tree[\graph] \triangleq \begin{cases}
      \dot{x}_\tree &= -\lap^{\tree} R_\graph W_\graph R_\graph^\intercal  x_\tree + D_\tree^\intercal E_{\graph}^{-1} \omega\\& - \lap^\tree R_\graph {v}\\
      z &= R_\graph^\intercal  x_\tree
                    \end{cases}\label{eq:2}\\
    \hat{\Sigma}^\tree[\graph] \triangleq \begin{cases}
      \dot{x}_\tree &= -\lap^{\tree} R_\graph W_\graph R_\graph^\intercal  x_\tree + D_\tree^\intercal E_{\graph}^{-1}\omega\\& - \lap^\tree R_\graph {v}\\
      z &=  x_\tree.
                    \end{cases}\label{eq:7}
\end{align}
There are multiple motivations for each of the above models, and the edge consensus formulation in general. 
As we can see in~\eqref{eq:16}, the edge state dynamics naturally exclude the zero eigenvalue of $\lap_\graph$.  
Also, as proposed in~\cite{Zelazo2011a}, it suffices to measure the $n-1$ relative states across the edges of a spanning tree $\tree$ to determine if the states of the $n$ nodes are at consensus.
Therefore, a satisfactory output measurement of the network is precisely the edge states of that spanning tree as in \eqref{eq:7}.
Of course, the remaining edges that are not in $\tree$ also influence the dynamics, and to study the $\mathcal{H}_2$ norm of the entire system one can take the output of these edge states into account as well.
The model in \eqref{eq:2} has output $z = R_\graph^T x_\tree$, which reconstructs the states of the remaining edges in $\graph\setminus\tree$ using the Tucker representation matrix, and thus considers all edges in the output.

In the next section, we discuss the $\Htwo$ norm as a system-theoretic metric of performance of these models.
For brevity, we may drop the argument of $\Sigma^\tree[\graph]$.

\subsection{Systems Theory and $\Htwo$ Performance}
\label{sec:systems-theory-htwo}
The input-output excitation properties of a linear system $\dot x = Ax + Bu,~y=Cx$ with transfer function $G(s) = C(sI-A)^{-1}B$ can be described using the $\Htwo$ norm,
\begin{align}
  \mathcal{H}_2^2(G) = \dfrac{1}{2\pi} \int_{-\infty}^\infty \mathbf{tr} \left[ G(j\omega)^{*}  G(j\omega) \right] d\omega.
\end{align}
The $\Htwo$ norm measures the steady-state covariance of the output of the system under zero-mean, unit-covariance white noise inputs, or equivalently the root-mean-square of the impulse response of the system.

As discussed in~ \cite{Foight2019a}, the $\Htwo$ norms of the models $\Sigma^\tree[\graph],\hat{\Sigma}^\tree[\graph]$ in \eqref{eq:2} and \eqref{eq:7} are obtained by solving the Lyapunov equation
\begin{align}
  &-\lap^\tree R_\graph W_\graph R_\graph^\intercal  X - X R_\graph W_\graph R_\graph^\intercal  \lap^\tree \\&+ D_\tree^\intercal  E_\graph^{-1}\Omega\Omega^\intercal  E_\graph^{-1} D_\tree + \lap^\tree R_\graph \Gamma \Gamma^\intercal  R_\graph^\intercal  \lap^\tree = 0.\label{eq:8}
\end{align}
An analytic solution to Equation~\eqref{eq:8} can be obtained with a convenient choice of the error covariances $\Omega$ and $\Gamma$, specifically $ \Omega = \sigma_\omega E_\graph^{1/2}$ and $\Gamma = \omega_v W_\graph^{1/2}$.
In this case, the solution of \eqref{eq:8} is given by
\begin{align}
  X^* = \dfrac{1}{2} \left [ \sigma_\omega^2 (R_\graph W_\graph R_\graph^\intercal  )^{-1} + \sigma_v^2 \lap^\tree \right],
\end{align}
and so  the $\Htwo$ norms of  $\Sigma^\tree[\graph],\hat{\Sigma}^\tree[\graph]$ are, respectively,
\begin{align}
 & \mathcal{H}_2^2(\Sigma^\tree[\graph]) = \mathbf{tr}\left[ R_\graph^\intercal  X^* R_\graph\right] \\ 
                        &= \dfrac{\sigma_\omega^2}{2} \mathbf{tr} \left[ R_\graph^\intercal  (R_\graph W_\graph R_\graph^\intercal  )^{-1} R_\graph \right] + \dfrac{\sigma_v^2}{2} \mathbf{tr} \left[ R_\graph^\intercal  \lap^\tree R_\graph \right] \label{eq:9}\\ 
 & \mathcal{H}_2^2(\hat{\Sigma}^\tree(\graph)) = \mathbf{tr} \left[ X^* \right]\\
& = \dfrac{\sigma_\omega^2}{2} \mathbf{tr} \left[  (R_\graph W_\graph R_\graph^\intercal  )^{-1}  \right] + \dfrac{\sigma_v^2}{2} \mathbf{tr} \left[  \lap^\tree\right].\label{eq:10}
\end{align}
Of course, the choice of the error covariances causes a loss of generality in the applicability of the analysis of this model.
In \cite{Foight2019a,Foight2019ts}, it was shown that this choice, plus the true error covariances can be combined with an ordering in the semi-definite cone of covariances to reasonably bound the true performance, suggesting that this model is a convenient proxy for general noise models. 

It was noted in \cite{Foight2019a} that the contributions of the weights and time scales can be separated in \eqref{eq:9} \& \eqref{eq:10} as, 
\begin{align}
  \Htwosq(\Sigma^\tree[\graph]) &= \Htwosq(\Sigma^\tree; W ) + \Htwosq(\Sigma^\tree; E )\\
  \Htwosq(\hat{\Sigma}^\tree[\graph]) &= \Htwosq(\hat{\Sigma}^\tree; W ) + \Htwosq(\hat{\Sigma}^\tree; E )
\end{align}
where
\begin{align}
  &\Htwosq(\Sigma^\tree; W) \triangleq  \dfrac{\sigma_\omega^2}{2} \mathbf{tr} \left[ R_\graph^\intercal  (R_\graph W_\graph R_\graph^\intercal  )^{-1} R_\graph \right]\\
  &\Htwosq(\hat{\Sigma}^\tree; W ) \triangleq \dfrac{\sigma_w^2}{2} \mathbf{tr} \left[  (R_\graph W_\graph R_\graph^\intercal  )^{-1}  \right] \\
  &\Htwosq(\Sigma^\tree; E ) \triangleq \dfrac{\sigma_v^2}{2} \mathbf{tr} \left[ R_\graph^\intercal  \lap^\tree R_\graph \right]\\
  &\Htwosq(\hat{\Sigma}^\tree; E ) \triangleq \dfrac{\sigma_v^2}{2} \mathbf{tr} \left[  \lap^\tree\right].
\end{align}
Furthermore, it was shown in \cite{Foight2019a} that the $\Htwo$ norm of $\Sigma$ can be written in terms of the $\Htwo$ norm of $\hat{\Sigma}$ as follows:
\begin{align}
  &\Htwosq(\Sigma^\tree; W) \\&= \Htwosq(\hat{\Sigma}^\tree; W) +  \dfrac{\sigma_w^2}{2} \mathrm{tr}\left[ T_\tree^\intercal (R_\graph W_\graph R_\graph)^{-1} T_\tree \right]\label{eq:17}\\
  &\Htwosq(\Sigma^\tree; E ) = \Htwosq(\hat{\Sigma}^\tree; E ) + \dfrac{\sigma_v^2}{2} \mathbf{tr}\left[T_\tree^\intercal L^\tree T_\tree \right].
\end{align}
We examine the behaviour of each of these components when an edge is added to a spanning tree -- $\Htwosq(\hat{\Sigma},W)$ and $\Htwosq(\hat{\Sigma},E)$ are discussed in \S\ref{sec:weight-cycle-select}, and $\Htwosq({\Sigma},E)$ \& $\Htwosq({\Sigma},W)$ in \S\ref{sec:contr-time-scal} and \S\ref{sec:contr-edges-mathc} respectively.


\section{Main Results}
\label{sec:main-results}

\subsection{Minimum-$\Htwo$ Spanning Tree}
\label{sec:minim-spann-tree-1}
In this section, we show how a minimum-$\mathcal{H}_2$ norm spanning tree of a graph $\graph$ can be found with a greedy algorithm.
Let $\graph$ be a weighted, time scaled graph with edge consensus dynamics~\eqref{eq:7} and corresponding $\Htwo$ norm \eqref{eq:10}.

Suppose that we want to remove communication links from $\graph$ until we are left with a minimally-connected graph --- a tree.
Such an operation can be used during a `stealth mode' for networked systems, during which time having the minimum number of communication links while still being connected is desired.
For example, a swarm of UAVs may be using consensus to agree on a formation heading, but may want to limit their chance of detection when entering a hostile area.
In this section, we show that one can choose the spanning tree $\tree\subseteq\graph$ that minimizes the $\Htwo$ norm in \eqref{eq:10} using a greedy algorithm.
This allows the network to minimize the number of communication links used in the consensus algorithm, while maintaining optimal noise rejection properties.
The main result of this section is summarized in Proposition~\ref{prop:mintree}.

\begin{algorithm}\label{alg.1}
\caption{Minimum--$\mathcal{H}_2$ Spanning Tree}\label{alg.discac}
\begin{algorithmic}[1]
\BState \emph{Input:} $\graph = (\mathcal{V},\mathcal{E},\mathcal{W},\mathcal{S})$ \newline
\BState \emph{Initialize:} 
\State Choose $v\in\mathcal{V}$, set $\mathcal{Q} = \mathcal{V} \setminus \{v\},\mathcal{P} = \{v\}$, $\mathcal{R} = \emptyset$
\While{ $|\mathcal{Q}| > 0$}
\State \emph{Set}
\begin{align}
  N(\mathcal{P}) = \{ j \in \mathcal{Q} ~:~i \in \mathcal{P},~ij\in\mathcal{E} \}.
\end{align}
\State \emph{Solve} 
\begin{align}
(\text{P}1):= \min_{j \in N(\mathcal{P})} \left\{ \sigma_v^2(\epsilon_k^{-1} + \epsilon_j^{-1}) + \sigma_\omega^2 w_{jk}^{-1} ~:~ jk \in \mathcal{E}\right\}
\end{align}
\State \emph{From ${\arg\min}$~}$(\text{P}1) := (j^*, j^*k)$, \emph{update}
\begin{align}
  \mathcal{P} &\mapsto \mathcal{P} \cup \{j^*\}\\
  \mathcal{R} &\mapsto \mathcal{R} \cup \{j^*k\}\\
  \mathcal{Q} &\mapsto \mathcal{Q}\setminus \{j^*\}
\end{align}
\EndWhile
\State \emph{Output:} $\tree^* = (\mathcal{P}, \mathcal{R}, \mathcal{W}(\mathcal{R}),\mathcal{S})$
\end{algorithmic}
\end{algorithm}

\begin{prop}
  \label{prop:mintree}
  Consider a graph $\graph$ with the edge consensus dynamics~\eqref{eq:7}, and with corresponding $\Htwo$ norm \eqref{eq:10}.
  Any spanning tree $\tree'\subseteq \graph$ has $\Htwo$ norm
  \begin{align}
    \mathcal{H}_2^2(\hat{\Sigma}^{\tree'}) = \dfrac{\sigma_\omega^2}{2} \sum_{e\in\tree'} w_e^{-1} + \dfrac{\sigma_v^2}{2} \sum_{i=1}^n \dfrac{\mathrm{deg}(i)}{\epsilon_i}.\label{eq:14}
  \end{align}

  Algorithm \ref{alg.discac} returns a spanning tree $\tree^*\subseteq \graph$ that minimizes \eqref{eq:14} out of all possible spanning trees $\tree\subseteq \graph$.
\end{prop}
\begin{pf}
  Define an auxiliary graph $\graph^\prime = (\Verts^\prime,\Edges^\prime,\Weights^\prime)$ on the same vertex and edge set as $\graph = (\Verts,\Edges,\Weights,\mathcal{S})$, (i.e., $\Verts^\prime = \Verts,~\Edges^\prime = \Edges$) but with edge weights $w(\graph^\prime) \in \Weights^\prime$ defined by
  \begin{align}
    w(\graph^\prime)_{ij} = \frac{\sigma_v^2}{\epsilon_i} + \frac{\sigma_v^2}{\epsilon_j} + \frac{\sigma_\omega^2}{w_{ij}},~\forall~i,j\in\Verts,~ij\in\Edges.\label{eq:3}
  \end{align}
  Then, define the \emph{total cost} of the graph $\graph^\prime$ as the sum of the edge weights:
  \begin{align}
    TC(\graph^\prime)  = \dfrac{1}{2}\sum_{e\in\Edges^\prime} w(\graph^\prime)_e. \label{eq:6}
  \end{align}
  It is clear from \eqref{eq:3} that
  \begin{align}
    TC(\graph^\prime) &= \dfrac{1}{2} \sum_{e\in\Edges^\prime} w(\graph^\prime)_e = \dfrac{1}{2}  \sum_{ij\in \Edges} \left[ \frac{\sigma_v^2}{\epsilon_i} + \frac{\sigma_v^2}{\epsilon_j} + \frac{\sigma_\omega^2}{w_{ij}} \right].\label{eq:4}
  \end{align}
  For every edge $ij$ we get one term with $\epsilon_i^{-1}$ in the sum of \eqref{eq:4}, and so summing over all edges yields $\mathrm{deg}(i)$ terms $\epsilon_i^{-1}$. 
  Hence, we can conclude that from \eqref{eq:10}, for any spanning tree $\tree'\subseteq \graph'$, its total cost is exactly the $\Htwo$ norm of the corresponding tree $\tree$ in $\graph$.
  \begin{align}
    TC(\tree^\prime) &= \dfrac{1}{2} \left[ \sum_{e\in\Edges} \frac{\sigma_\omega^2}{w_e} + \sum_{i\in \Verts} \dfrac{\sigma_v^2\mathrm{deg}(i)}{\epsilon_i}\right]= \mathcal{H}_2^2(\Sigma^\tree).\label{eq:5}
  \end{align}
  By contruction, the minimal spanning tree $\tree^{\prime*}$ of $\graph^\prime$, defined by
  \begin{align}
    \tree^{\prime*} := \arg\min \left\{ TC(\tree^\prime) ~:~\tree^\prime \text{ is a spanning tree of }\graph^\prime  \right\},
  \end{align}
  corresponds precisely to the spanning tree $\tree$ of $\graph$ minimizing the $\Htwo$ norm in \eqref{eq:10}.

  It is well-known that a minimal-weighted spanning tree of $\graph^\prime$ (i.e., one that minimizes \eqref{eq:6}) can be obtained by the Jarn\'{i}k-Dijkstra-Prim (JDP), or Kruskal algorithms; see \cite{algorithms2011}.
  Algorithm~\ref{alg.discac} is precisely a JDP algorithm that builds the minimal-weighted spanning tree $\tree^*$ by starting with a single node in $\tree^*$, and at each iteration adding a new node to $\tree^*$ by picking the neighbouring vertex $k$ to $\tree^*$ (but $k$ not already in $\mathcal{T}^*)$ with the smallest term $\sigma_\omega^2w_{jk}^{-1} + \sigma_v^2(\epsilon_j^{-1} + \epsilon_k^{-1})$ for $j\in\tree^*$.
  By the above argument, this JDP algorithm on $\graph'$ with cost $ TC(\graph^\prime) $ in \eqref{eq:6} outputs the minimum-$\Htwo$ norm spanning tree of $\graph$.
 \qed
\end{pf}

\begin{rem}
  For an edge-weighted graph where the cost is the sum of the edge weights, it is known that if the edge weights are distinct, then the minimal-weighted spanning tree is unique.
  The minimum-$\Htwo$ spanning tree $\tree^*$ generated by Algorithm~\ref{alg.discac} is not unique, even if all the nodal time scales and edge weights are distinct.
  See \S\ref{sec:minim-spann-tree} for an example.
\end{rem}

\begin{rem}
Another practical algorithm would be a distributed version of Algorithm~\ref{alg.discac}, which could be constructed by considering a distributed minimal-weighted spanning tree algorithm such as \cite{Awerbuch1987,Gallager1983,Kutten1998}.
Such an algorithm would allow a consensus network to autonomously restructure itself to the $\Htwo$-optimal tree configuration in a distributed manner.
\end{rem}

\subsection{Weighted Cycle Selection}
\label{sec:weight-cycle-select}
In the previous section, we discussed how one may find the minimum-$\mathcal{H}_2$ norm spanning tree from a connected graph $\graph$.
In this section, we examine what happens to the $\mathcal{H}_2$ norm when edges are added to a spanning tree.
In particular, we show that adding weighted edges improve the $\mathcal{H}_2$ norm, and discuss what choice of edge optimizes this improvement.
This can be used in the previous motivating example of a stealthy system that wants to still minimize the number of communication links, but does not reject noise well enough with just the links in the spanning tree.

\cite{Zelazo2013} discussed how adding cycles to an unweighted, mono-scaled graph impacts the $\Htwo$ performance in the presence of unit covariance noise applied to the edges and nodes.
In particular, they found that the biggest improvement in $\Htwo$ performance resulted from adding an edge to a tree $\tree$ which maximized the length of the resulting cycle. 
Our first contribution is showing that in the presence of edge weights and time scales, the length of the cycle no longer matters, rather one should add an edge to form a cycle that, roughly speaking, has small weights.

First, we define the \emph{weighted length} and \emph{unweighted length} of a cycle $\cycle$ as,
\begin{align}
  l_w(\cycle) := \sum_{e\in \cycle} w_e^{-1},~l(\cycle) := \sum_{e\in \cycle} 1.
\end{align}

\begin{prop}
  \label{prop:edge1}
  Consider a graph consisting of a weighted, time scaled tree $\tree$, and consider the dynamics $\hat{\Sigma}^\tree$, with the $\Htwo$ norm of $\hat{\Sigma}^\tree$ given by Equation~\eqref{eq:10}.
  Consider the task of adding edges to $\tree$ to optimize the $\Htwo$ norm of $\hat{\Sigma}^\tree$.
  Then, the following hold:
\begin{enumerate}
  \item Adding an edge $e$ with weight $W_{\cycle}$ to $\tree$ improves the $\Htwo$ performance in the following manner: 
\begin{align}
 \mathcal{H}_2^2(\hat{\Sigma}[\tree\cup e])=   \mathcal{H}_2^2(\hat{\Sigma}[\tree]) - \frac{\sigma_\omega^2}{2l_w(\cycle)}\sum_{ij\in\tree\cap\cycle} w_{ij}^{-2},~~~\label{eq:20}
\end{align}
where $\cycle$ is the unique cycle formed by adding edge $e$ to $\tree$.
Furthermore, adding an edge to $\tree$ to form a cycle always decreases the $\Htwo$ norm.
  \item Adding $p$ edge-disjoint cycles $\cycle_1,\dots,\cycle_p$ via edges $e_1,\dots,e_p$ improves the $\Htwo$ performance in the following manner:
    \begin{align}
      & \mathcal{H}_2^2(\hat{\Sigma}[\tree\cup \{e_i\}])\\&=   \mathcal{H}_2^2(\hat{\Sigma}[\tree]) - \sum_{k=1}^p \left[ \frac{\sigma_\omega^2}{2l_w(\cycle_k)} \sum_{ij\in\tree\cap\cycle_k} w_{ij}^{-2} \right],
    \end{align}
\end{enumerate}
\end{prop}
\begin{pf}
  Recall that the $\Htwo$ norm of \eqref{eq:7} for a given spanning tree $\tree\subseteq \graph$ is given by
  \begin{align}
    \mathcal{H}_2^2(\hat{\Sigma}^\tree[\graph]) = \dfrac{\sigma_\omega^2}{2} \mathbf{tr}\left[ ( R_\graph W_\graph R_\graph^\intercal  )^{-1} \right] + \dfrac{\sigma_v^2}{2} \mathbf{tr} \left[ \lap^\tree\right].
  \end{align}
  Since the graph in question is an edge added to a tree, we have $\graph = \tree \cup e$.
  This introduces one cycle, $\cycle$.
  Our choice of spanning tree in the cycle representation matrix $T_\tree$ is the original tree $\tree$ to which the edge $e$ is added.
  Consider the trace argument in the weight matrix term, 
  \begin{align}
    &( R_\graph W_\graph R_\graph^\intercal  )^{-1}  =  \left(
    \begin{bmatrix}
      I & T_\tree 
    \end{bmatrix}
    \begin{bmatrix}
      W_\tree &  \mathbf{0} \\ \mathbf{0} & W_\cycle
    \end{bmatrix}
    \begin{bmatrix}
      I \\ T_\tree^\intercal 
    \end{bmatrix}\right)^{-1}\\
                                           &=    \left(W_\tree + T_\tree  W_\cycle T_\tree^\intercal  \right)^{-1}
                                           = \left(W_\tree + T^\prime T^{\prime \intercal}\right)^{-1},
  \end{align}
  where $W_\tree$ is the top-left block of $W_\graph$ corresponding to the weights of the edges in $\tree$, and $W_\cycle$ is the bottom-right block of $W_\graph$ corresponding to the weight of the edge being added to $\tree$.
  In the last display, we made the substitution $T^\prime \triangleq T_\tree  W_\cycle^{1/2}$.

Since we are adding a single edge $e$ to $\tree$,  $W_\cycle$ is a scalar and $T^\prime$ is a vector.
Therefore, $T^\prime T^{\prime \intercal}$ is a rank-one matrix, and we can apply the Sherman-Morrison-Woodbury update, yielding
\begin{align}
  ( R_\graph W_\graph R_\graph^\intercal  )^{-1}  &= W_\tree^{-1} - \dfrac{W_\tree^{-1} T^\prime T^{\prime \intercal} W_\tree^{-1}}{1+ T^{\prime \intercal} W_\tree^{-1} T^\prime}\\
                                         &= W_\tree^{-1} - \dfrac{W_\tree^{-1} T_\tree  W_\cycle T_\tree ^\intercal  W_\tree^{-1}}{1 + W_\cycle \left(T_\tree^\intercal W_\tree^{-1} T_\tree \right)}.\label{eq:11}
\end{align}
We recall the following lemma.
\begin{lem}[Prop 1 \& 2 in \cite{Zelazo2013}]
  \label{lem:zelazo1}
  Recall the definition of the cycle representation matrix:
  \begin{align}
    T_\tree = (D_\tree^\intercal D_\tree)^{-1} D_\tree^\intercal D_{\graph \setminus \tree} := \left[
    \begin{matrix}
      c_1 & & \cdots & & c_{|\Edges(\graph\setminus\tree)|}
    \end{matrix}\right].
  \end{align}
  Each column $c_i$ of $T_\tree$ represents a cycle $\cycle_i$ of $\graph$.
  The matrices $T_\tree^\intercal T_\tree$ and $T_\tree T_\tree^\intercal$ encode the following information about the cycles of $\graph$:
  \begin{enumerate}
    \item $\left[T_\tree^\intercal T_\tree\right]_{ii} = l(\cycle_i) - 1$
    \item $\left[T_\tree^\intercal T_\tree\right]_{ij} = 0$ if and only if the cycles $c_i$ and $c_j$ are edge-disjoint.
    \item $\left[T_\tree T_\tree^\intercal\right]_{ee}$ is the number of times edge $e$ is used to construct the cycles of $\graph$.
  \end{enumerate}
\end{lem}

Consider the graph in question, $G:= \tree \cup e$.
Adding this edge introduces a single cycle, and so by Lemma~\ref{lem:zelazo1},
\begin{align}
  T^\intercal_\tree T_\tree = c_1^Tc_1 = l(\cycle) - 1.
\end{align}
This holds because $T^\intercal_\tree T_\tree$ adds a `1' to the sum for each edge in the tree that ends up in the cycle $\mathcal{C}$, in other words we can write
\begin{align}
  T^\intercal_\tree T_\tree = \left[\sum_{e\in \tree\cap\cycle} (1)\right]  = \left[\sum_{e\in\cycle} (1)\right] - 1,
\end{align}
where the `-1' comes from excluding the edge $e$ that is added to $\tree$.
Similarly, it can be seen that 
\begin{align}
     &(T_\tree )^\intercal  W_\tree^{-1}T_\tree  = \sum_{ij\in\cycle \cap \tree} w_{ij}^{-1}\\
     & = \left[\sum_{ij\in\cycle} w_{ij}^{-1}\right] - W_\cycle^{-1} = l_w(\cycle) - W_\cycle^{-1}.
\end{align}
From this, we can conclude that
\begin{align}
  &1 + W_\cycle (T_\tree )^\intercal W_\tree^{-1} T_\tree  = 1 + W_\cycle T^\intercal W_\tree^{-1}T\\
                                                             &= 1 + W_\cycle \left( l_w(\cycle) - W_\cycle ^{-1} \right)
  = W_\cycle l_w(\cycle),
\end{align}
and so from \eqref{eq:11} we have that
\begin{align}
     ( R_\graph W_\graph R_\graph^\intercal  )^{-1}  &=  W_\tree^{-1} - \dfrac{1}{l_w(\cycle)} W_\tree^{-1} TT^\intercal  W_\tree^{-1},\label{eq:13}
\end{align}
where the length of the cycle added by edge $e$ is given by
\begin{align}
    l_w(\cycle) = \sum_{e\in\tree\cap\cycle} w_e^{-1} + W_\cycle^{-1}.
\end{align}
Recall from  Lemma~\ref{lem:zelazo1} that $[T_\tree T_\tree^\intercal ]_{ee}$ is the number of times the edge $e\in\tree$ is used in the cycle $\cycle$.
Therefore, 
\begin{align}
  \left[W_\tree^{-1} T_\tree T_\tree^\intercal  W_\tree^{-1}\right]_{ee} = w_e^{-2}.\label{eq:12}
\end{align}

Taking the trace of \eqref{eq:13} and using \eqref{eq:12} yields
\begin{align}
  & \mathbf{tr}\left[( R_\graph W_\graph R_\graph^\intercal  )^{-1}\right]\\
  &= \mathbf{tr}\left[ W_\tree^{-1} \right] + \dfrac{1}{l_w(\cycle)}\mathbf{tr}\left[W_\tree^{-1} TT^\intercal  W_\tree^{-1}\right]\\
                           &= \frac{1}{l_w(\cycle)} \left[ \sum_{e\in\tree\cap\cycle} w_e^{-2} \right].
\end{align}
For the second part, note that the improvement in $\Htwo$ performance after adding edge $c$ is
\begin{align}
 - \dfrac{\sigma_\omega^2}{2l_w(\cycle)} \left[ \sum_{ e\in \tree\cap\cycle}  w_e^{-2} \right] < 0,
\end{align}
since $l_w(\cycle) >0 $ and $W_e^{-2} >0$ for all $e\in\Edges$, and so the $\Htwo$ norm always decreases.

In the case of adding $k$-edge disjoint cycles, each edge can only be in one cycle, and so by Lemma~\ref{lem:zelazo1} we have that \eqref{eq:12} still holds.
The result follows by applying the preceding argument for each of the $k$ cycles.
\qed
\end{pf}

\subsection{Contributions of Time Scales to $\mathcal{H}_2^2(\Sigma^\tree[\tree])$}
\label{sec:contr-time-scal}
As discussed in \S\ref{sec:systems-theory-htwo}, the two models $\Sigma^\tree$, $\hat{\Sigma}^\tree$ differ in the outputs --  the former considers only the edge states of the spanning tree, and the latter considers all  edges of the graph.
We noted in the proof of Proposition~\ref{prop:mintree} that when adding an edge to a tree, the $\mathcal{H}_2$ norm of the dynamics $\hat{\Sigma}^\tree$ is not affected by the time scales of the nodes on which the edge is placed.
In this section, we show that these dynamics , adding an edge to slow nodes provides the smallest increase in $\Htwo$ norm.

Recall that the $\Htwo$ performance of $\Sigma^\tree$ can be written as,
\begin{align}
  \mathcal{H}_2^2(\Sigma^\tree) &= \frac{\sigma_\omega^2}{2} \mathbf{tr} [R_\graph^\intercal  (R_\graph W_\graph R_\graph^\intercal )^{-1} R] + \frac{\sigma_v^2}{2} \mathbf{tr}[R_\graph^\intercal  \lap^\tree R_\graph]\\
                &:= \mathcal{H}_2^2(\Sigma^\tree;W) + \mathcal{H}_2^2(\Sigma^\tree; E).\label{eq:15}
\end{align}
Here we have separated the contributions of the weights and time scales into two terms.
Consider the second term in \eqref{eq:15}, defined as 
\begin{align}
  \mathcal{H}_2^2(\Sigma^\tree; E) := \frac{1}{2} \sigma_v^2 \mathbf{tr}[R_\graph^\intercal  \lap^\tree R_\graph].
\end{align}
We have the following proposition.
\begin{prop}
  \label{prop:ts}
  There exists a separation of $\Htwo(\Sigma, E)$ into two terms, one depending on the cycle states, and the other the spanning tree states.
  Namely,
  \begin{align}
    \mathcal{H}_2^2(\Sigma; E) =\dfrac{\sigma^2_v}{2}\left( \Tr[\lap^\tree] + \mathbf{tr}[\lap^{\graph\setminus\tree}]\right).
  \end{align}
\end{prop}
\begin{pf}
  By trace cyclicity, we can compute: 
  \begin{align}
    &\mathcal{H}_2^2\mathrm(\Sigma^\tree, E) = \Tr[R_\graph^\intercal \lap^\tree R_\graph]\\
                     &= \Tr\left[
                       \begin{bmatrix}
                         I \\ (T_\tree)^\intercal 
                       \end{bmatrix}
                       \lap^\tree
    \begin{bmatrix}
      I & T_\tree
    \end{bmatrix}
                       \right]
                     = \Tr[\lap^\tree] + \Tr[T_\tree^\intercal \lap^\tree T_\tree].
  \end{align}
  Via trace cyclicity, we can manipulate the argument of the second trace term in the last display as follows:
  \begin{align}
    &\Tr[T_\tree^\intercal \lap^\tree T_\tree]= \\
    & \Tr[D_\tree (D_\tree ^\intercal D_\tree )^{-1}D_\tree ^\intercal D_{\graph\setminus\tree} D_{\graph\setminus\tree}^\intercal D_\tree  (D_\tree ^\intercal D_\tree )^{-1}D_\tree ^\intercal E_\graph^{-1}  ] 
  \end{align}
  Notice that the term $D_\tree  (D_\tree ^\intercal D_\tree )^{-1} D_\tree ^\intercal $ is the projection operator on the range space of $D_\tree $.
  Since the cycles are linear combinations of the spanning tree edges, every column of $D_{\graph\setminus\tree}$ is in the range of $D_\tree $.
  Therefore, the projection  of $D_{\graph\setminus\tree}$ onto the range space of $D_\tree $ is nothing but $D_{\graph\setminus\tree}$:
  \begin{align}
    D_\tree  (D_\tree ^\intercal D_\tree )^{-1} D_\tree ^\intercal  D_{\graph\setminus\tree} = D_{\graph\setminus\tree},
  \end{align}
  and so we can conclude that
  \begin{align}
    &\Tr[T_\tree^\intercal \lap^\tree T_\tree ] = 
                                                 \Tr[D_{\graph\setminus\tree} D_{\graph\setminus\tree}^\intercal D_\tree  (D_\tree ^\intercal D_\tree )^{-1}D_\tree ^\intercal E_\graph^{-1}  ] \\
   &= \Tr[D_{\graph\setminus\tree}^\intercal E_\graph^{-1} D_{\graph\setminus\tree} ] = \Tr[\lap^{\graph\setminus\tree}].
  \end{align}
\qed
\end{pf}

From Proposition~\ref{prop:ts}, we can conclude with the following corollary.
\begin{cor}
  Consider the dynamics~\eqref{eq:2} of $\Sigma^\tree$ on a graph $\graph$.
  Then, consider adding an edge $ij$ to $\graph$ to get $\graph\cup\{ij\}$.
  The time scale term of the $\mathcal{H}_2^2$ norm in~\eqref{eq:15} satisfies
  \begin{align}
      \mathcal{H}_2^2(\Sigma^\tree[\graph \cup \{ij\}];E )  = \dfrac{\sigma_v^2}{2} \left(\epsilon_i^{-1} + \epsilon_j^{-1} \right) + \mathcal{H}_2^2 (\Sigma^\tree[\graph],E ).
  \end{align}
\end{cor}

\subsection{Contribution of Weights to $\mathcal{H}_2^2(\Sigma^\tree[\tree])$}
\label{sec:contr-edges-mathc}
In \S\ref{sec:weight-cycle-select} we showed that adding an edge to a spanning tree always improves the $\mathcal{H}_2$ norm of $\hat{\Sigma}$.
In this section, we discuss what happens to $\mathcal{H}_2$ norm of $\Sigma$ by examining \eqref{eq:17}.

Recall from \S\ref{sec:systems-theory-htwo} that we can write
\begin{align}
    \Htwosq(\Sigma^\tree; W) &= \Htwosq(\hat{\Sigma}^\tree; W)\\& +  \dfrac{\sigma_w^2}{2} \mathrm{tr}\left[ T_\tree^\intercal (R_\graph W_\graph R_\graph)^{-1} T_\tree \right].\label{eq:18}
\end{align}
When $\graph = \tree$, we examined what happens to $\mathcal{H}_2^2(\hat{\Sigma}^\tree; W )$ when an edge is added to $\tree$ in Proposition~\ref{prop:edge1}.
When $\graph = \tree$, the second term in~\eqref{eq:18} is zero since $T_\tree$ is the empty matrix.
In the next proposition, we show that this term is positive when an edge is added to $\tree$.

\begin{prop}
  Consider the edge consensus model \eqref{eq:2} on a graph $\tree$ with $\Htwo$ norm \eqref{eq:9}.
  When an edge $e$ with weight $W_\cycle$ is added to $\tree$ forming a cycle $\cycle$, the weight portion of the $\Htwo$ norm satisfies
  \begin{align}
      \Htwosq(\Sigma^\tree[\tree\cup\{ij\}];W) &= \Htwosq(\Sigma^\tree[\tree];W)\\ &+\frac{\sigma_\omega^2}{2 W_\cycle}  - \frac{\sigma_\omega^2}{2l_w(\cycle)}  \sum_{ij\in\cycle} w_{ij}^{-2}. \label{eq:21}
  \end{align}
  \begin{pf}
    First, note that the term in \eqref{eq:18} is zero when there are no cycles in $\graph$.
    Adding a single edge to the tree $\tree$ forms a cycle, and hence $T_\tree$ has one column.
    We can compute using the rank-one update in \eqref{eq:11}, 
    \begin{align}
      &\dfrac{\sigma_\omega^2}{2} \mathbf{tr}\left[ T_\tree^\intercal (R_\graph W_\graph R_\graph)^{-1} T_\tree \right]\\
      &= \dfrac{\sigma_\omega^2}{2} \mathbf{tr}\left[T_\tree^\intercal \left(W_\tree^{-1} - \frac{1}{l_W(\cycle)} W_\tree^{-1} T_\tree T_\tree^\intercal W_\tree^{-1}\right)T_\tree\right]\\
      &= \dfrac{\sigma_\omega^2}{2} \left(\mathbf{tr} \left[T_\tree ^\intercal W_\tree^{-1} T_\tree \right]  - \frac{1}{l_w(\cycle)} \mathbf{tr}\left[(T_\tree ^\intercal W_\tree^{-1} T_\tree)^2\right] \right)\\
      &= \dfrac{\sigma_\omega}{2} \left(l_w(\cycle) -1 - \frac{(l_w(\cycle)-1)^2}{l_w(\cycle)} \right)\\
      &= \dfrac{\sigma_\omega}{2}\left( \dfrac{l_w(\cycle) W_\cycle^{-1} - W_\cycle^{-2}}{l_w(\cycle)} \right).\label{eq:19}
    \end{align}
      Summing \eqref{eq:19} with the improvement of $\Htwosq(\hat{\Sigma};W)$ in \eqref{eq:20} from Proposition~\ref{prop:edge1} yields~\eqref{eq:21}.
      \qed
    \end{pf}

\end{prop}


\section{Examples}
\label{sec:examples}

\subsection{Minimum $\Htwo$-Norm Spanning Trees}
\label{sec:minim-spann-tree}

\begin{figure}
    \centering
    \begin{subfigure}[b]{\columnwidth}
      \vspace{-2mm}
      \includegraphics[width=0.9\columnwidth]{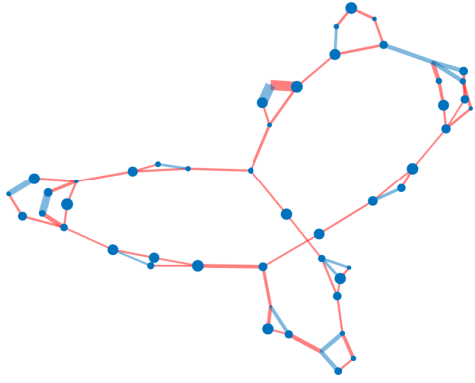}
      \caption{Series-parallel network}
      \label{fig:sp}
    \end{subfigure}
    \begin{subfigure}[b]{\columnwidth}
      \includegraphics[width=0.9\columnwidth]{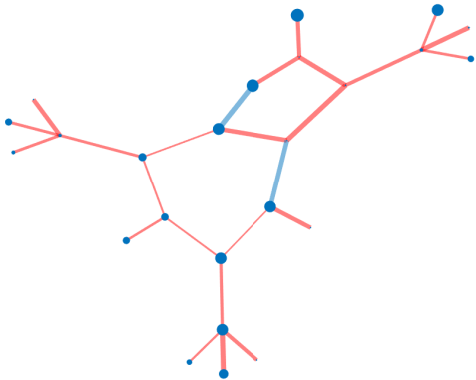}
      \caption{Adjacency graph of the caffeine molecule}
      \label{fig:caf}
    \end{subfigure}
    \begin{subfigure}[b]{\columnwidth}
      \includegraphics[width=0.9\columnwidth]{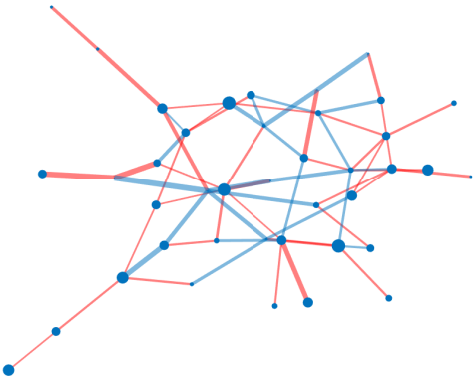}
      \caption{Erd\H{o}s-R\'{e}nyi random graph with $P(ij\in\mathcal{E}) = 0.09$}
      \label{fig:rand}
    \end{subfigure}
    \caption{Example graphs with the minimum-$\Htwo$ norm spanning tree (highlighted in red), computed from Algorithm~\ref{alg.discac}. Rejected edges from the original graph are highlighted in blue. Edge widths are proportional to $\log(w_e)$, and node markers are proportional to $\epsilon_i$.}
    \label{fig:tre4s}
  \end{figure}

In this section, we show examples of the minimum-$\Htwo$ spanning trees from Algorithm~\ref{alg.discac}.
Three such graphs are shown in Figure~\ref{fig:tre4s}; edges in the minimum-$\Htwo$ spanning tree are highlighted in red, and the remaining rejected edges in blue.
In each example, the edge weights $w_e$ and time scales $\epsilon_i$ are drawn from a uniform random distribution on $[0,1]$ --- the width of each edge is proportional to $\log (w_e)$ for clarity, and the node size is proportional to $\epsilon_i$.
The graph in Figure~\ref{fig:sp} is an example of a series-parallel network; the graph in Figure~\ref{fig:caf} is the molecular adjacency structure of the caffeine molecule, and the graph in Figure~\ref{fig:rand} is an Erd\H{o}s-R\'{e}nyi random graph.
The random graph is generated on 40 nodes, and an edge between node $i$ and $j$ exists with probability 0.09.

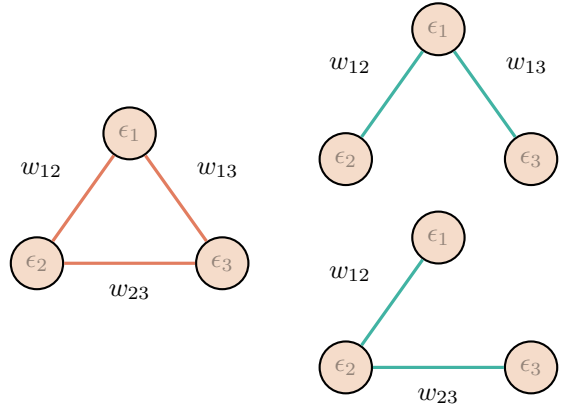
\begin{figure}
  \centering
  \begin{scaletikzpicturetowidth}{0.8\columnwidth}
    \begin{tikzpicture}[scale=\tikzscale]
      \begin{scope}[every node/.style={circle,thick,draw,fill=nodecl,fill opacity=0.4,text opacity = 1}]
        \node (2) at (0,0) {$\epsilon_2$};
        \node (1) at (1.5,2.12) {$\epsilon_1$};
        \node (3) at (3,0) {$\epsilon_3$};
        \node (22) at (5,1.7) {$\epsilon_2$};
        \node (11) at (6.5,3.82) {$\epsilon_1$};
        \node (33) at (8,1.7) {$\epsilon_3$};
        \node (222) at (5,-1.7) {$\epsilon_2$};
        \node (111) at (6.5,0.42) {$\epsilon_1$};
        \node (333) at (8,-1.7) {$\epsilon_3$};
      \end{scope}

      \begin{scope}[>={stealth[black]},
        every edge/.style={draw=edgecl,very thick}]
        \path [-] (1) edge node [label=above left:{$w_{12}$}] {} (2);
        \path [-] (1) edge node [label=above right:{$w_{13}$}] {} (3); 
        \path [-] (2) edge node [label=below :{$w_{23}$}] {} (3);
      \end{scope}

      \begin{scope}[>={stealth[black]},
        every edge/.style={draw=treecl,very thick}]
        \path [-] (11) edge node [label=above left:{$w_{12}$}] {} (22);
        \path [-] (11) edge node [label=above right:{$w_{13}$}] {} (33); 
        \path [-] (111) edge node [label=above left:{$w_{12}$}] {} (222);
        \path [-] (222) edge node [label=below :{$w_{23}$}] {} (333); 
      \end{scope}

    \end{tikzpicture}
  \end{scaletikzpicturetowidth}
  \caption{Example of a time scaled, weighted graph with distinct edge weights and time scales, but no unique minimal $\Htwo$ spanning tree.}
  \label{fig:nonunique}
\end{figure}

Next, we show an example of a graph with distinct time scales on the nodes, and distinct edge weights, but no unique minimum-$\mathcal{H}_2$ spanning tree.
Consider the triangle in Figure~\ref{fig:nonunique} with the edge consensus model $\hat{\Sigma}^\tree$, and the following time scales and edge weights:
\begin{align}
  \epsilon_1 &=1,~w_{23} = 1\\
  \epsilon_2 &=2,~w_{13} = 2\\
  \epsilon_3 &=3,~w_{12} = 3.
\end{align}
Further suppose that the noise weightings satisfy $\sigma_\omega = \sigma_v = 1$.
Then, the spanning tree $\tree_1$ of the triangle in Figure~\ref{fig:nonunique} on the top right has $\Htwo$ norm
\begin{align}
  \Htwo(\Sigma^{\tree_1} (\graph)) & = \dfrac{1}{2} \left[ w_{12}^{-1} + w_{13}^{-1} +  \epsilon_1^{-1} + \sum_{i=1}^3 \epsilon_i^2 \right] = 1.8\bar{3}
\end{align}
and the spanning tree $\tree_2$  in Figure~\ref{fig:nonunique} on the bottom right has $\Htwo$ norm
\begin{align}
  \Htwo(\Sigma^{\tree_2} (\graph)) & = \dfrac{1}{2} \left[ w_{12}^{-1} + w_{23}^{-1} +  \epsilon_3^{-1} + \sum_{i=1}^3 \epsilon_i^2 \right]= 1.8\bar{3}.
\end{align}
\begin{rem}
  The reason the minimum-$\mathcal{H}_2$ spanning trees in Figure~\ref{fig:nonunique} are not unique, despite having distinct edge weights and time scales in $\graph$, is because the edge weights of the auxilary graph $\graph'$ defined in the proof of Proposition~\ref{prop:mintree}, given by
  \begin{align}
    w^\prime_{ij} = \sigma_v^2(\epsilon_i^{-1} + \epsilon_j^{-1}) + \sigma_\omega^2w_{ij}^{-1},
  \end{align}
  are not distinct in $\graph'$.
\end{rem}

\subsection{Adding Edges to a Tree}

\begin{figure}
    \centering
    \begin{scaletikzpicturetowidth}{0.9\columnwidth}
      \begin{tikzpicture}[scale=\tikzscale]
        \begin{scope}[every node/.style={circle,thick,draw,fill=nodecl,fill opacity=0.4,text opacity = 1}]
          \node (2) at (0,0) {};
          \node (1) at (1.5,2.12) {};
          \node (3) at (3,0) {};
          \node (4) at (4.5,2.12) {};
          \node (5) at (7,2.12) {};
          \node (6) at (8.5,0) {};
        \end{scope}

        \begin{scope}[>={stealth[black]},
          every edge/.style={draw=edgecl,very thick}]
          \path [-] (1) edge  (2);
          \path [-] (1) edge  (3); 
          \path [-] (3) edge  (4); 
          \path [-] (4) edge  (5); 
          \path [-] (5) edge  (6); 
        \end{scope}

        \begin{scope}[>={stealth[black]},
          every edge/.style={draw=treecl,very thick}]
          \path [-] (2) edge node [label=below :{$W_1$}] {} (3); 
          \path [-] (3) edge node [label=below :{$W_{2}$}] {} (6);
        \end{scope}

      \end{tikzpicture}
    \end{scaletikzpicturetowidth}
    \caption{Example of a weighted graph where a cycle with fewer edges has a lower $\mathcal{H}_2$ norm for $\hat{\Sigma}^\tree$ than a cycle with more edges.}
    \label{fig:cyclechoice}
\end{figure}
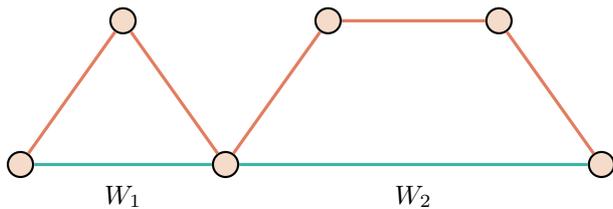

Consider the path graph in Figure~\ref{fig:cyclechoice} denoted by the solid edges, where all the edges in the path have weight $w_e = 1$.
Suppose $W_1 = 10$ and $W_2=5$, and $\sigma_\omega =1$
Then, the improvement of adding either edge $W_1$ to $\graph$ to create cycle $\cycle_1$, or adding $W_2$ to $\graph$ to create cycle $\cycle_2$ is
\begin{align}
     - \dfrac{1}{2l_w(\cycle_1)} \left[ \sum_{ e\in \tree\cap\cycle_1}  w_e^{-2} \right] &\approx -0.9524\\
     - \dfrac{1}{2l_w(\cycle_2)} \left[ \sum_{ e\in \tree\cap\cycle_2}  w_e^{-2} \right] &\approx -0.9375,
\end{align}
and so the shorter (but more highly weighted) cycle $\cycle_1$ offers the better improvement.

\section{Conclusion}
\label{sec:conclusion}

In this paper, we have considered various perspectives on optimizing weights and time scales for edge consensus.
We showed that a greedy algorithm can be used to find the minimum-$\Htwo$ norm spanning tree of a graph, and then discussed how the $\Htwo$ norm can be optimized by adding cycles for two output models of edge consensus.

Future work includes implementing a distributed version of Algorithm~\ref{alg.discac}, as well as examining how to optimally choose more general graph structures.
Choosing time scales to improve controllability or observability of edge consensus is also an open problem.


\section*{ACKNOWLEDGEMENTS}
This  work  was  supported  by  the  Swiss  Competence  Centers  for  Energy  Research  FEEB\&D project  and  the  ETH  Foundation. 

\bibliography{ConsensusKoopman}

\end{document}